\documentclass{lms}

\usepackage[notref,notcite,color,final]{showkeys}
\usepackage{amsmath}
\usepackage{amssymb}
\usepackage[utf8]{inputenc}

\newcommand{\N}{\mathbb{N}} 
 
\newcommand{\R}{\mathbb{R}} 
\newcommand{\C}{\mathbb{C}} 
\newcommand{\Schwartz}{\mathcal{S}}
\newcommand{\Fspace}{\Schwartz}
\newcommand{\omegaclass}{\overline{\omega}}
\newcommand{\paren}[1]{{\left( #1 \right)}} 
 
\newcommand{\abs}[1]{{\left\lvert #1 \right\rvert}} 
\newcommand{\bigabs}[1]{{\bigl\lvert#1\bigr\rvert}} %
\newcommand{\norm}[2]{{\left\lVert #2 \right\rVert}_{#1}} 
\newcommand{\bignorm}[2]{{\bigl\lVert #2 \bigr\rVert}_{#1}} 
\newcommand{\norma}[2]{{\left\lvert #2 \right\rvert}_{#1}} 
\newcommand{\normb}[2]{{\left\lvert\!\left\lvert\!\left\lvert #2 
\right\rvert\!\right\rvert\!\right\rvert}_{#1}} 
\newcommand{\scal}[2]{\left\langle #1, #2 \right\rangle} 
\newcommand{\deq}{\mathrel{\mathop:}=} 
\newcommand{\Fourier}[1]{\widehat{#1}}
\newcommand{\suchthat}{, } 
\newcommand{\set}[2]{\left\{#1 \suchthat #2 \right\}}
\newcommand{\prob}{\mathbf{p}}
\newcommand{\floor}[1]{\left\lfloor#1\right\rfloor}
\newcommand{\ceil}[1]{\left\lceil#1\right\rceil}
\DeclareMathOperator{\trace}{tr}
\newtheorem{theo}{Theorem} 
\newtheorem{prop}{Proposition}[section] 
\newtheorem{coro}[prop]{Corollary} 
\newtheorem{lemm}[prop]{Lemma} 
\newnumbered{defi}{Definition}

\title[Ultrarapidly decreasing ultradifferentiable
functions]{Ultrarapidly decreasing ultradifferentiable functions,
Wigner distributions and density matrices}

\author{Jean-Marie Aubry}

\classno{Primary 46A61; Secondary 33C45, 46A13.}

\extraline{The author acknowledges the support of the French Agence
Nationale de la Recherche (ANR), under grant StatQuant (JCJC07 205763)
``Quantum Statistics''.}

\simpleequations

\begin{document}

\maketitle

\begin{abstract}
    Spaces $\Fspace_{\omega}, \Fspace_{\{\omegaclass\}},
    \Fspace_{(\omegaclass)}$ of ultradecreasing ultradifferentiable
    (or for short, ultra-$\Schwartz$) functions, depending on a weight
    $e^{\omega(x)}$, are introduced in the context of quantum
    statistics.  The corresponding coefficient spaces in the Fock
    basis are identified, and it is shown that the Hermite expansion
    is a tame isomorphism between these spaces.  These results are
    used to link decrease properties of density matrices to
    corresponding properties of the Wigner distribution.
\end{abstract}

\section{Introduction}

\subsection{Quantum states}

A quantum state is a vector of unit norm in some Hilbert space, for
instance $L^2(\R)$.  In quantum optics (our subjacent model throughout
this paper), the simplest states decrease like the Gaussian, and have
the same regularity (their Fourier transform decreases like the
Gaussian as well).  One can however easily encounter \emph{super
Gaussian} states whose wave function decreases still rapidly, but not
as fast as the Gaussian, like for instance $\abs{f(x)} \leq C
e^{-\frac{\abs{x}^\beta}{2}}$ for $0 < \beta \leq 2$, and whose
Fourier transform decreases at the same rate (the restriction $\beta
\leq 2$ comes from the uncertainty principle, detailed
in~\S~\ref{sec:uncertaintyprinciple}).  Rapidly decreasing
self-Fourier functions, of importance in
optics~\cite{Caola:1991qy,Corcoran:2004ul,Dragoman:1996mz}, fall under
this category.

The class of compactly supported functions $f$ that have a regularity
comparable to $\bigabs{\Fourier{f}(\xi)} \leq C
e^{-\frac{\abs{\xi}^\beta}{2}}$ is called a class of
\emph{ultradifferentiable} functions (in the sense of Beurling,
see~\cite{Braun:1990rt}).  The classes of functions that we introduce
in~\S~\ref{sec:ultraspace} are related but strictly larger, because we
replace compact support with an ``ultrafast'' decrease condition.

When both conditions ($\abs{f(x)} \leq C e^{-\frac{\abs{x}^\beta}{2}}$
and $\bigabs{\Fourier{f}(\xi)} \leq C e^{-\frac{\abs{\xi}^\beta}{2}}$)
are satisfied, $f$ can be called an ``ultrafast decreasing
ultradifferentiable'' function.  More generally, following 
Gr\"ochenig and Zimmerman~\cite{Grochenig:2004fk}, we can replace $x
\mapsto \frac{\abs{x}^\beta}{2}$ by a \emph{weight function} $\omega$
satisfying certain conditions (Definition~\ref{defi:weightfunction}).
The spaces thus obtained can be viewed as the ultradifferentiable
version of the Schwartz class $\Schwartz$; for short we call them
spaces of \emph{ultra-$\Schwartz$} functions.

\subsection{The Fock basis}

The natural basis of $L^2(\R)$ for the problems of quantum optics is
the one formed by the Hermite functions
\begin{equation*}
    h_m(x) \deq \paren{2^m m!  \sqrt\pi}^{-\frac12} H_m(x)
    e^{-\frac{x^2}{2}},
\end{equation*}
where $H_m(x) \deq (-1)^m e^{x^2} \frac{d^m}{dx^m} e^{-x^2}$ is the
$m$-th Hermite polynomial.  The normalization above ensures that
$\scal{h_{m}}{h_{n}} = \delta_{m,n}$.

In terms of time/frequency, the Hermite functions are perfectly
balanced, because they are eigenvectors of the Fourier
transform\footnote{The definition that we use is $\Fourier{f}(\xi)
\deq (2\pi)^{-\frac{1}{2}} \int f(x) e^{-i\xi x} dx$}
\begin{equation}
    \label{eq:FourierHermite}
    \Fourier{h_{m}} = (-i)^m h_{m}.
\end{equation}
This basis is also called the \emph{Fock basis} and the $h_{m}$ are
the \emph{Fock states} (the $m$-photons state in the case of light). 
Note that (\ref{eq:FourierHermite}) characterizes ``self-Fourier'' 
functions as the ones having non-zero Hermite coefficients only for 
$m \in 4\N$.

A few more properties will be needed. For instance, $h_{n}$ is also an 
eigenvector of the \emph{Hermite operator}
\begin{equation}
    (x^2-\tfrac{d^2}{dx^2}) h_{n}(x) = (2n+1) h_{n}(x).
    \nonumber
\end{equation}
Some decrease properties of $h_{n}(x)$ are stated in 
Lemma~\ref{lemm:dh}.

\subsection{Outline of the results}

Despite their simple definition, the study of the functional
properties of Hermite series is in many cases a difficult problem.
The $L^p$ case, for instance, is the object of a book by
Thangavelu~\cite{Thangavelu:1993fk}.  A question similar to ours has
also been investigated numerically by Boyd in~\cite{Boyd:1984uq}.  The
case $\omega(x) = C x$ was treated by Janssen and van
Eijndhoven~\cite{Janssen:1990fk}.  More anecdotically, the decrease of
the Hermite coefficients of a certain function has been related to
properties of the Riemann $\zeta$ function by Grawe
in~\cite{Grawe:1997uq}.  Finally,
Langenbruch~\cite{Langenbruch:2006ly} has treated the problem starting
from a slightly different definition for the function spaces
(conditions of the type $\norm{\infty}{x^\alpha D^{\beta} f} \leq
M_{\alpha+\beta}$); in fact a part of his proofs could be adapted and
reused here.

Our first goal in this paper is to characterize the class of
ultra-$\Schwartz$ functions introduced above in terms of Hermite
coefficients; this is done, up to a multiplicative constant on
$\omega$, in~\S~\ref{sec:Hermiteexpansion}.  The main result in this
paper is Theorem~\ref{theo:mainHermite} in \S~\ref{sec:mainHermite}.

The next question occurs when one considers a so-called {mixed state}:
that is, the statistical result of a physical experiment where there
is an uncertainty on the quantum state itself.  Let us say that, with
probabilities $\prob_{0},\prob_{1},\dots$, the experiment produces
orthogonal states $\psi_{0},\psi_{1},\dots$.  All this information is
contained in the operator
\begin{equation*}
    \rho \deq \sum_{i} \prob_{i} \scal{.}{\psi_{i}} \psi_{i},
\end{equation*}
where $\scal{.}{\psi_{i}} \psi_{i}$ is the projector on the subspace
generated by $\psi_{i}$.  If more than one of the $\prob_{i}$ is
non-zero, then the state is called \emph{mixed}.  In other words, a
pure state is one that can be described by a single wave functions:
for instance $\psi = \sum_{i} \sqrt{\prob_{i}} \psi_{i}$.

It is easily seen, as in Leonhardt's reference
book~\cite{Leonhardt:1997fk}, that $\rho$ is Hermitian, that its
eigenvalues are nonnegative with a sum $\trace(\rho) = 1$, and that
these conditions are sufficient to make a mixed quantum state.

For a more graphic representation, one can also consider the
\emph{Wigner distribution} associated to $\rho$.  It is a function of
two variables $q,p \in \R$ (or rather, a function of $q+ip \in \C$),
defined as
\begin{equation}
    \label{eq:Wignerdistribution}
    \Phi_{\rho}(q,p) \deq \sum_{i} \frac{\prob_{i}}{2\pi} \int e^{ipx} 
    \psi_{i}(q-\frac{x}{2}) \overline{\psi_{i}(q+\frac{x}{2})} dx.
\end{equation}
For some reasons, it is sometimes easier to work with 
\begin{equation*}
    \tilde\Phi_{\rho} (q,p) \deq \frac12
    \Phi_{\rho}\paren{\frac{q}{\sqrt2},\frac{p}{\sqrt2}}.
\end{equation*}
The transform $\rho \mapsto \tilde\Phi_{\rho}$ is an isometry between
$\mathcal{L}(L^2(\R))$ and $L^2(\C)$; its inverse is called the
\emph{Weyl transform}\footnote{See~\cite{Thangavelu:1993fk} for more
insight into this deep connexion with the Heisenberg group
representations.}.  The Wigner distribution is very useful for the
interpretation of experiments, such as in quantum homodyne tomography
({\it op.  cit.}, see also the lectures notes edited by Paris and
\v{R}eh\'a\v{c}ek~\cite{Paris:2004fk}).

For a pure state $\psi$, or in terms of operators, $\rho_{\psi} =
\scal{.}{\psi}\psi$, the Wigner representation of $\rho_{\psi}$ is
given by a simpler form of~(\ref{eq:Wignerdistribution}) called the
\emph{Wigner transform} of $\psi$, which is the quadratic form
\begin{equation}
    \label{eq:Wignertransform}
    \Phi(\psi,\psi)(q,p) \deq \Phi_{\rho_{\psi}}(q,p) = \frac{1}{2\pi}
    \int e^{ipx} \psi(q-\frac{x}{2}) \overline{\psi(q+\frac{x}{2})} dx.
\end{equation}
Then there is no difficulty to see that $\int \Phi(\psi,\psi)(q,p) dp
= \abs{\psi(q)}^2$ and that $\int \Phi(\psi,\psi)(q,p) dq =
\bigabs{\Fourier\psi(p)}^2$; this is what makes this representation
useful in quantum mechanics, for it allows to recover the probability
distribution of an observable and its dual ({\it e.g.} positition and
momentum, electric and magnetic field\ldots) as marginal distributions
of the Wigner transform of $\psi$ (even though $\Phi(\psi,\psi)$,
being in general non-positive, is not a probability density
\textit{stricto sensu}).

Another consequence is that, generally speaking, a decrease condition
on $\Phi(\psi,\psi)$ implies a decrease condition on both $\psi$ and
$\Fourier{\psi}$.  From the results of~\S~\ref{sec:Hermiteexpansion},
this implies a decrease condition on the Hermite coefficients of
$\psi$: the object of~\S~\ref{sec:Wd} is to show that the
converse is true.

This result can be extended to mixed states in the following manner:
let $[\rho_{m,n}]$ be the (infinite) \emph{density matrix} of $\rho$
in the Fock basis.  For instance, in the case of a pure state $\rho = 
\rho_{\psi}$,
\begin{equation}
    \rho_{m,n} = \scal{\psi}{h_{m}} \overline{\scal{\psi}{h_{n}}}.
    \label{eq:rhomnforpurestates}
\end{equation}
Assuming that $\rho_{m,n}$ decreases ultrarapidly when $m+n\to\infty$,
we shall show that $\Phi_{\rho}$ also decreases ultrarapidly when
$\abs{p}+\abs{q}\to\infty$.  However, in \S~\ref{sec:mixedstates} we
see that there are (necessarily mixed) states for which $\Phi_{\rho}$
decreases ultrarapidly but not $\rho_{m,n}$.

Another object of interest is $\Fourier{\Phi_{\rho}}$, the
2-dimensional Fourier transform of $\Phi_{\rho}$.  In radar
technology, this function is known as the \emph{ambiguity function}.
Its decrease rate rules the regularity of $\Phi_{\rho}$, which is of
importance in statistical estimation, see Butucea {\it et
al.}~\cite{Butucea:2005uq}.  We study the decrease and regularity
properties of $\Fourier{\Phi_{\rho}}$ in
\S~\ref{sec:ambiguityfunction}.

This leads to studying the matrix elements of the Weyl transform:
combining~(\ref{eq:Wignerdistribution}) and~(\ref{eq:Wignertransform})
we see that, writing
\begin{equation}
    \label{eq:defPhimn}
    \Phi_{m,n}(q,p) \deq \Phi(h_{m},h_{n})(q,p) = \frac{1}{2\pi} \int e^{ipx}
    h_{m}(q-\frac{x}{2}) {h_{n}(q+\frac{x}{2})} dx,
\end{equation}
we have just $\Phi_{\rho} = \sum_{m,n} \rho_{m,n}\Phi_{m,n}$.  Similarly,
we define
\begin{equation*}
    \tilde\Phi_{m,n} (q,p) \deq \frac12
    \Phi_{m,n}\paren{\frac{q}{\sqrt2},\frac{p}{\sqrt2}}.
\end{equation*}
These functions, named \emph{special Hermite functions} by Strichartz
~\cite{Thangavelu:1993fk}, are remarkable.  Although they are not
exactly 2-dimensional Hermite functions (the latter being defined by
tensor product $h_{u, v}(q, p)\deq h_{u}(q) h_{v}(p)$), special
Hermite functions are a special linear combination (for $u+v = m+n$)
of those.  Consequently, they are also eigenfunctions for the
2-dimensional Fourier transform:
\begin{equation}
    \Fourier{\tilde\Phi_{m,n}} = (-i)^{m+n} 
    \tilde\Phi_{m,n}.
    \label{eq:FourierspecialHermite}
\end{equation}

Let us now expose our frame of work.

\section{Functional setting}
\label{sec:ultraspace}

\subsection{Ultra-$\Schwartz$ spaces}

Because of (\ref{eq:FourierHermite}), any function space for which the
Fock states (Hermite functions) form an unconditional basis must be
invariant under Fourier transform.  Such spaces can be defined by
weighted $L^\infty$ conditions on both $f$ and $\Fourier{f}$, as
follows.

\begin{defi}
    \label{defi:weightfunction}
    A continuous increasing function $\omega: [0,\infty) \to
    [0,\infty)$ is called a \emph{weight function} if it satisfies
    \begin{enumerate}
	\renewcommand{\theenumi}{(\roman{enumi})}
	\renewcommand{\labelenumi}{\theenumi}
	
	\item $\Omega : t \mapsto \omega(e^t)$ is convex
	\label{it:convexOmega}
	\label{it:firstomega}
	
	\item $\log(t) = o(\omega(t))$ 
	\label{it:minomega}
	
	\item $\omega(2t) = O(\omega(t))$ \label{it:croissomega}
	
	\item $\limsup \frac{\omega(t)}{t^2} \leq \frac{1}{2}$.
	\label{it:maxomega}
	\label{it:lastomega}

	\newcounter{enumaxomega}
	\setcounter{enumaxomega}{\value{enumi}}
	
    \end{enumerate}
\end{defi}

This definition should be compared with that of Braun, Meise and
Taylor~\cite{Braun:1990rt}.  The only difference is that their
($\beta$) $\int_{1}^{\infty}\frac{\omega(t)}{t^{2}}dt < \infty$ is
replaced by~\ref{it:maxomega}, which is weaker; this is because we
relax the compact support condition.  As we shall see
below,~\ref{it:maxomega} is related to the uncertainty principle; for
technical reasons, we shall sometimes need instead the strict
inequality {\it
   \begin{enumerate}
       \setcounter{enumi}{\value{enumaxomega}}
       \renewcommand{\theenumi}{(\roman{enumi})}
       \renewcommand{\labelenumi}{\theenumi} 
       
       \item $\limsup \frac{\omega(t)}{t^2} < \frac{1}{2}$.
       \label{it:maxomegastrict}
   \end{enumerate}}

Given a weight function $\omega$ satisfying
\ref{it:firstomega}--\ref{it:lastomega}, we define
\begin{equation*}
    \Fspace_{\omega} \deq \set{f \in L^2(\R)}{ f(x) = O(e^{-
    \omega(\abs{x})}) \text{ and } \Fourier{f}(\xi) = O(e^{-
    \omega(\abs{\xi})})}.
\end{equation*}
Clearly $\Fspace_{\omega}$ is a Banach space, if equipped with the
norm
\begin{equation*}
    \norm{\omega}{f} \deq \norm{\infty}{f(x) e^{\omega(\abs{x})}} +
    \norm{\infty}{\Fourier{f}(\xi) e^{\omega(\abs{\xi})} }.
\end{equation*}

Remark that $\lambda \mapsto \Fspace_{\lambda\omega}$ is strictly
decreasing (in the sense of inclusion) on $(0,+\infty)$.
Nevertheless, it will be convenient to state our results in terms of
spaces that depend only on the rate of growth of $\omega$, that is, on
the class $\omegaclass \deq \set{\lambda\omega}{\lambda > 0}$.  For
this purpose, we introduce the ultra-$\Fspace$ spaces
\begin{equation*}
    \Fspace_{\{\omegaclass\}} \deq \bigcup_{\epsilon > 0}
    \Fspace_{\epsilon \omega}
\end{equation*}
and
\begin{equation*}
    \Fspace_{(\omegaclass)} \deq \bigcap_{N<\infty} \Fspace_{N\omega}.
\end{equation*}
Note that the constant in~\ref{it:croissomega} depends on
$\omegaclass$ only, we shall often use this fact.  Also remark that in
the above union and intersection, the parameter can be taken in a
countable set without changing anything to the definition.

\subsection{Uncertainty principle}
\label{sec:uncertaintyprinciple}

Bearing the name of Heisenberg, the principle that says that one
cannot know with precision both the position and the momentum of a
particle translates quantitatively in terms of simultaneous
localization of a function and its Fourier transform; for a
comprehensive study follow Havin and J\"oricke~\cite{Havin:1997fk}.
In our setting, this principle says that if $\omega$ increases too
fast, then $\Fspace_{\omega}$ is trivial (reduced to $\{0\}$).

More precisely, let us quote the one-dimensional version of a theorem
of Bonami, Demange and Jaming \cite{Bonami:2003fk}, in the descent of
Hardy~\cite{Hardy:1933lr}, Beurling and
H\"ormander~\cite{Hormander:1991fk}.
\begin{theo*}[(BDJ)]
    Let $f \in L^2(\R)$ and $N \geq 0$. Then
    \begin{equation}
	\iint\limits_{\R\times\R} \frac{\abs{f(x)}
	\bigabs{\Fourier{f}(\xi)} e^{\abs{x \xi}} }
	{\paren{1+\abs{x}+\abs{\xi}}^N} dx d\xi < \infty \nonumber
    \end{equation}
    if and only if $f(x) = p(x) e^{-\frac{x^{2}}{2}}$, where $p$ is a
    polynom of degree $< \frac{N-1}{2}$.
\end{theo*}

This theorem means that if
\begin{equation}
    N_{\omega} \deq \inf \set{N}{ \iint\limits_{\R^+\times\R^+}
    \frac{e^{-\paren{\omega(x) + \omega(\xi) -x \xi } } }
    {\paren{1+{x}+{\xi}}^N} dx d\xi < \infty } \nonumber
\end{equation}
then the dimension of $S_{\omega}$ is 
\begin{equation}
\ceil{\frac{N_{\omega} - 1}{2}} \leq \dim(S_{\omega}) \leq
\floor{\frac{N_{\omega}+1}{2}}.  \nonumber
\end{equation}
When $\frac{N_{\omega}+1}{2} \in \N$, then one has to check whether
$\iint\limits_{\R^+\times\R^+} \frac{e^{-\paren{\omega(x) +
\omega(\xi) -x \xi } } } {\paren{1+{x}+{\xi}}^{N_{\omega}}} dx d\xi$
is finite or not, to see whether $\dim(S_{\omega}) =
\frac{N_{\omega}-1}{2}$ or $\frac{N_{\omega}+1}{2}$ respectively.  In
particular, the necessary and sufficient condition for $S_{\omega}$ to
be non trivial is that
\begin{equation}
    \iint\limits_{\R^+\times\R^+} 
    e^{-\paren{\omega(x)+\omega(\xi)-x\xi}} dx d\xi = \infty.
    \nonumber
\end{equation}

The ``critical case'' $\omega(x)=\frac{x^{2}}{2}$ is however a fuzzy
frontier, because it is actually the limit inferior of
$\frac{\omega(x)}{x^2}$ that counts when determining whether the above
integral diverges or not.  Indeed it is easy to construct a weight
function $\omega$ such that $\liminf \frac{\omega(x)}{x^2} <
\frac{1}{2} < \limsup \frac{\omega(x)}{x^2}$ (the only thing to take
care of is~\ref{it:convexOmega}, a piecewise affine function $\Omega$
with increasing slopes will do).  This shows that~\ref{it:maxomega}, a
fortiori~\ref{it:maxomegastrict}, is indeed a restriction to the range
of our results.\footnote{In fact~\ref{it:maxomega} is only used in
Lemma \ref{lemm:numueff} to guarantee that the constant $K$ is
universal, but any fixed number could replace $\frac12$.  Theorem
\ref{theo:mainHermite} still holds if we replace~\ref{it:maxomega} by
a finite limit superior.  However, when~\ref{it:maxomegastrict} is
asked for in Propositions~\ref{prop:lb},
\ref{prop:Wl},~\ref{prop:Wlmixed}, then $\frac{1}{2}$ is indeed the
critical value.}
There might be some room for improvement here.

Another consequence of Theorem (BDJ) is that a sufficient condition
for $ \Fspace_{\{\omegaclass\}}$ to be non-trivial is given by
$\limsup \frac{\omega(t)}{t^2} < \infty$; for
$\Fspace_{(\omegaclass)}$ it is that $\limsup \frac{\omega(t)}{t^2} =
0$.

For comparison, we recall that the fastest possible decrease for the
Fourier transform of compactly supported functions is of order $e^{- C
\abs{t}}$, see~\cite{Aubry:tf,Braun:1990rt} for instance.

\subsection{Sequence spaces}

The spaces of Hermite coefficients for functions in
$\Fspace_{\{\omegaclass\}}$ and $\Fspace_{(\omegaclass)}$ will be
identified in~\S~\ref{sec:Hermiteexpansion} to the following spaces.
Let $\omega$ be a weight function as in
Definition~\ref{defi:weightfunction}.  Then
\begin{equation}
    \Lambda_{\omega} \deq \set{(\alpha_{n}) \in\C^\N}{\alpha_{n} =
    O(e^{-\omega(\sqrt{n})})} \nonumber
\end{equation}
endowed with the norm $\norm{\omega}{(\alpha_{n})} \deq \sup_{n}
\abs{\alpha_{n}} e^{\omega(\sqrt n)}$ is a Banach space, and as
previously we define
\begin{equation}
    \Lambda_{\{\omegaclass\}} \deq \bigcup_{\epsilon>0}
    \Lambda_{\epsilon\omega} \nonumber
\end{equation}
as well as
\begin{equation}
    \Lambda_{(\omegaclass)} \deq \bigcap_{N<\infty} \Lambda_{N\omega}.
    \nonumber
\end{equation}

\subsection{Topology}

As a decreasing intersection of Banach spaces,
$\Fspace_{(\omegaclass)}$ endowed with the projective topology (the
coarsest topology that makes every embedding $\Fspace_{(\omegaclass)}
\to \Fspace_{N\omega}$ continuous) is a Fr\'echet space.  The same
holds for $\Lambda_{(\omegaclass)}$.

The picture is a little more complicated for
$\Fspace_{\{\omegaclass\}}$ and $\Lambda_{\{\omegaclass\}}$.  As an
increasing union of Banach spaces, ($\Lambda_{\alpha \omega} \subset
\Lambda_{\beta\omega}$ if $\beta < \alpha$),
$\Lambda_{\{\omegaclass\}}$ is naturally endowed with the inductive
limit topology (the finest locally convex topology such that every
embedding $\Lambda_{\epsilon \omega} \to \Lambda_{\{\omegaclass\}}$ is
continuous).  In the standard terminology, $\Lambda_{\{\omegaclass\}}$
is called a (LB)-space.  The same goes for
$\Fspace_{\{\omegaclass\}}$.  It is a classical problem in functional
analysis to study the properties of such spaces, in particular their
completeness.  Note that this inductive limit is not strict (in the
sense of K\"othe~\cite[\S~19.4]{Kothe:1969fr}) because the topology on
$\Lambda_{\alpha \omega}$ is not induced by that of
$\Lambda_{\beta\omega}$.

\begin{prop}
    \label{prop:Lambdacomplete}
    $\Lambda_{\{\omegaclass\}}$ is a Silva space, thus complete.
\end{prop}
    
\begin{proof}
    Recall, as in~\cite{Cooper:1969qy}, that a Silva space is a
    locally convex inductive countable union of increasing Banach
    spaces, such that each embedding map is compact (it is in
    particular a (DFS)-space).  This is the case for
    $\Lambda_{\{\omegaclass\}}$, as we can see that the union may be
    taken over $\epsilon \in \set{\frac1m}{m\in\N}$, that
    $\Lambda_{\frac{\omega}{m} } \subset \Lambda_{\frac{\omega}{m+1}
    }$ and that the unit ball of $\Lambda_{\frac{\omega}{m} }$ is
    relatively compact in $\Lambda_{\frac{\omega}{m+1} }$ because
    $e^{\frac{\omega(\sqrt{n})}{m}} = o\paren{e^{
    \frac{\omega(\sqrt{n})}{m+1} }}$ as $n \to \infty$.
    
    \nocite{Bierstedt:2003yq} \nocite{Jarchow:1981sp} \nocite{Meise:1997fk}
\end{proof}

\begin{coro}
    $\Fspace_{\{\omegaclass\}}$ is complete.
\end{coro}

\begin{proof}
    Anticipating a little, by Theorem~\ref{theo:mainHermite},
    $\Fspace_{\{\omegaclass\}}$ is the image of
    $\Lambda_{\{\omegaclass\}}$ by a topological isomorphism.
\end{proof}

\subsection{Tame isomorphisms}

The notion of \emph{tameness} was introduced
in~\cite{Hamilton:1982fk}, in connexion with the Nash-Moser inverse
function theorem.  It is a natural notion of regularity for linear
operators between Fr\'echet spaces or inductive limits of Banach
spaces.

\newcommand{\Cclass}{\overline{C}}

\begin{defi}
    \label{defi:tame}
    Let $(E_{j}, \abs{\ }_{j})_{j\in\N}$ and $(F_{j}, \norm{j}{\
    })_{j\in\N}$ be two increasing families of Banach spaces, and let
    $E\deq \bigcup_{j}E_{j}$, $F \deq \bigcup_{j}F_{j}$ be the
    corresponding (LB)-spaces endowed with the inductive topologies.
    A linear mapping $T: E \to F$ is called \emph{tame} if there are
    $j_{0}\in \N$, $\Cclass < \infty$ such that for all $j \geq
    j_{0}$, there exists $C < \infty$ satisfying
    \begin{equation}
	\norm{\Cclass j}{T(f)} \leq C \abs{f}_{j}.
	\label{eq:tameinductive}
    \end{equation}
    
    Let $(E_{j}, \abs{\ }_{j})_{j\in\N}$ and $(F_{j}, \norm{j}{\
    })_{j\in\N}$ be two decreasing families of Banach spaces (or more
    generally, spaces with increasing semi-norms), and let $E\deq
    \bigcap_{j}E_{j}$, $F \deq \bigcap_{j}F_{j}$ be the corresponding
    Fr\'echet spaces endowed with the projective topologies.  In that
    case, a linear mapping $T: E \to F$ is said to be \emph{tame} if
    there are $j_{0}\in \N$, $\Cclass < \infty$ such that for all $j
    \geq j_{0}$, there exists $C < \infty$ satisfying
    \begin{equation}
	\norm{j}{T(f)} \leq C \abs{f}_{\Cclass j}.
	\label{eq:tameprojective}
    \end{equation}
    
    A linear mapping $T$ is called a \emph{tame isomorphism} if it is
    bijective and if both $T$ and $T^{-1}$ are tame.
\end{defi}

\section{Hermite expansion}
\label{sec:Hermiteexpansion}

\subsection{Main results}

In our setting, (\ref{eq:tameinductive}) of Definition~\ref{defi:tame}
will apply to spaces of type $\{\omegaclass\}$
and~(\ref{eq:tameprojective}) will apply to spaces of type
$(\omegaclass)$.

\label{sec:mainHermite}
\begin{theo}
    \label{theo:mainHermite}
    Let $H: f \mapsto (\scal{f}{h_{n}})_{n\in\N}$.
    
    If $\limsup \frac{\omega(t)}{t^2} < \infty$, then $H$ is a tame
    isomorphism: $\Fspace_{\{\omegaclass\}} \to
    \Lambda_{\{\omegaclass\}}$.
    
   If $\limsup \frac{\omega(t)}{t^2} = 0$, then $H$ is a tame
   isomorphism: $\Fspace_{(\omegaclass)} \to \Lambda_{(\omegaclass)}$.
\end{theo}

This theorem is proved in two parts, first the upper bound:
\begin{prop}
    \label{prop:ub}
    There exist $\Cclass < \infty$, depending only on $\omegaclass$, and
    $C < \infty$ (which may depend on $\omega$) such that, for all
    $f \in \Fspace_{\Cclass \omega}$,
    \begin{equation}
	\norm{\omega}{H(f)} \leq
	C \norm{\Cclass \omega}{f}
	\label{eq:ub}
    \end{equation}
\end{prop}
and the lower bound:
\begin{prop}
    \label{prop:lb}
    Assume~\ref{it:maxomegastrict}.  There exist $\Cclass < \infty$,
    depending only on $\omegaclass$, and $C < \infty$ (which may
    depend on $\omega$) such that, for all $\alpha \in
    \Lambda_{\Cclass \omega}$,
    \begin{equation}
	\norm{ \omega }{H^{-1}(\alpha)} \leq C \norm{\Cclass
	\omega}{\alpha}.
	\label{eq:lb}
    \end{equation}
\end{prop}

The two inequalities above give immediately the projective
case~(\ref{eq:tameprojective}), with $\abs{f}_{j} \deq \norm{j
\omega}{f}$ and $\norm{j}{H(f) } \deq \norm{j \omega}{H(f) }$.  The
condition $\limsup \frac{\omega(t)}{t^2} = 0$ ensures that
Proposition~\ref{prop:lb} can be applied to any $j\omega$.

Replacing $\omega$ by ${\Cclass}^{-1} \omega$ in the previous two
propositions, we get~(\ref{eq:tameinductive}), with $\abs{f}_{j} \deq
\norm{\frac \omega j}{f}$ and $\norm{j}{H(f) } \deq \norm{\frac \omega
j}{H(f)}$.  In that case, the condition $\limsup \frac{\omega(t)}{t^2}
< \infty$ guarantees that
Proposition~\ref{prop:lb} can be applied to
$\frac{\omega}{j}$ when $j$ is large enough.  Thus
Theorem~\ref{theo:mainHermite} is proved.

We now turn to the proof of Propositions \ref{prop:ub} and
\ref{prop:lb}.  In the sequel, $C$ denotes a constant (not necessarily
persistent) which may depend on $\omega$ and $\Cclass$ denotes a
constant which may depend only on $\omegaclass$.  In Lemmas
\ref{lemm:numueff} and \ref{lemm:Hermiteoperator}, $K$ denotes a
universal constant.

\subsection{Upper bound}

The proof is based on a few elementary lemmata.  We consider that a
weight function $\omega$ has been fixed satisfying Definition
\ref{defi:weightfunction}, and we recall that $\omegaclass \deq
\set{\lambda\omega}{\lambda > 0}$.

\begin{lemm}
    \label{lemm:Omegastar}
    Let $\Omega : t \mapsto \omega(e^t)$ and let $\Omega^\star : \nu
    \mapsto \sup_{t\in\R} \nu t - \Omega(t)$ be its convex conjugate.
    There exists a $C < \infty$ such that, for all $f
    \in\Fspace_{\omega}$, for all $\nu \geq 0$,
    \begin{equation*}
	\norm{L^2}{x^\nu f} \leq C e^{\Omega^\star(\nu+1)}
	\norm{\omega}{f}
    \end{equation*}
    and    (here $D = \frac{d}{dx}$)
    \begin{equation*}
	\norm{L^2}{D^\mu f} \leq C e^{\Omega^\star(\mu+1)}
	\norm{\omega}{f}.
    \end{equation*}

\end{lemm}
Note that~\ref{it:minomega} ensures $\Omega^\star(\nu) < \infty$ for
all $\nu \geq 0$.

\begin{proof}
    Since $\norm{\omega}{f} = \bignorm{\omega}{\Fourier{f}}$, we only
    have to prove the first inequality.
    
    Note that $\Omega^\star$ is increasing, so
    \begin{align*}
	\int_{\abs x < 1} \abs{x}^{2\nu}\abs{f(x)}^2 dx &\leq 2
	\sup_{\abs x < 1} \abs{f(x)}^2 \leq 2 \norm{\omega}{f}^2 \\
	 &\leq 2 {e^{-2\Omega^\star(1)}} e^{2\Omega^\star(\nu+1)}
	 \norm{\omega}{f}^2.
    \end{align*}    
	
    On the other hand we have
    \begin{align*}
	\int_{\abs{x}>1} \abs{x}^{2\nu} \abs{f(x)}^2 dx & \leq \sup_{x
	> 1} e^{(2\nu+2)\log(x)-2\omega(x)} \norm{\omega}{f}^2
	\int_{\abs{x}>1} \frac{dx}{x^2} \\
	&\leq 2 e^{2 \sup_{t\geq 0} (\nu+1) t -\Omega(t)}
	\norm{\omega}{f}^2 \\
	&\leq 2 e^{2 \Omega^\star(\nu+1)} \norm{\omega}{f}^2.
    \end{align*}
    Adding up the two, we get the result (with $C =
    \sqrt2(e^{-\Omega^\star(1)} + 1)$).
\end{proof}

The second lemma does the interpolation between the information on $f$
and the information on $\Fourier{f}$.

\begin{lemm}
    \label{lemm:numueff}
    There exists a $K < \infty$ such that, for all weight functions
    $\omega$, there exists $C$, for all $f \in \Fspace_{\omega}$, for
    all $\mu , \nu \in\N$,
    \begin{equation}
	\label{eq:numueff}
	\norm{L^2}{x^\nu D^\mu f} \leq C K^{\mu+\nu}
	e^{\frac{\Omega^\star(2\mu+2\nu+2)}{2}} \norm{\omega}{f}.
    \end{equation}
\end{lemm}

\begin{proof}
    We write
    \begin{align*}
	\norm{L^2}{x^\nu D^\mu f}^2 &= \int D^\mu \overline{f(x)}
	x^{2\nu} D^\mu f(x) dx \\
	&= (-1)^\mu \int \overline{f(x)} D^\mu \paren{x^{2\nu} D^\mu
	f(x)} dx \\
	&= (-1)^\mu \int \overline{f(x)} \sum_{k=0}^{\mu \wedge 2\nu}
	\binom{\mu}{k} \frac{(2\nu)!}{(2\nu-k)!} x^{2\nu-k} D^{2\mu-k}
	f(x) dx \\
	&\leq \sum_{k=0}^{\mu \wedge 2\nu} \binom{\mu}{k}
	\binom{2\nu}{k} k!  \norm{L^2}{x^{2\nu-k} f}
	\norm{L^2}{D^{2\mu-k} f} \\
	&\leq C \sum_{k=0}^{\mu \wedge 2\nu}
	\binom{\mu}{k} \binom{2\nu}{k} k!
	e^{\Omega^\star(2\nu-k+1)+\Omega^\star(2\mu-k+1)}
	\norm{\omega}{f}^2
    \end{align*}
    thanks to Lemma \ref{lemm:Omegastar} and the fact that
    $\Fspace_{\omega} \subset \Fspace$.
    
    By convexity, $\nu \mapsto \Omega^\star(\nu)-\Omega^\star(0)$ is
    superadditive, so $\Omega^\star(2\nu-k+1) \leq
    \Omega^\star(2\nu+1) - \Omega^\star(k) +\Omega^\star(0)$, and the
    same with $\mu$.  Thanks to~\ref{it:maxomega}, there exists $x_{0}
    =: e^{t_{0}}$ such that $x \geq x_{0} \Rightarrow \omega(x) \leq
    x^2$.  It follows that
    \begin{align*}
	\Omega^\star(k) &\geq \sup_{t < t_{0}} kt - \Omega(t) \vee
	\sup_{t \geq t_{0}} k t - e^{2t} \\
	&\geq \frac{k}{2} \paren{\log\paren{\frac{k}{2 }}-1}
    \end{align*}
    when $k \geq k_{0} \deq 2 x_{0}^{2}$.  In that case, using
    Stirling's formula, we see that there exists a constant $B$ such
    that $\frac{\log(k!)}{2} - \Omega^\star(k) \leq k \frac B 2$.
    
    At this point we have shown that when $k \geq k_{0}$,
    \begin{align*}
	\log(k!)  + \Omega^\star(2\nu-k+1) \quad & \\
	+ \Omega^\star(2\mu-k+1) & \leq \Omega^\star(2\nu+1) +
	\Omega^\star(2\mu+1) + kB + 2 \Omega^\star(0)  \\
	\intertext{using the superadditivity again,}
	&\leq \Omega^\star(2\nu+2\mu+2) + kB + 3 \Omega^\star(0)
    \end{align*}
    we crudely bound $\binom{\mu}{k} \binom{2\nu}{k}$ by $4^{\mu+\nu}$
    and $e^{kB}$ by
    $\paren{e^{2B}}^{\mu+\nu}$, and the lemma follows.
\end{proof}

To finish, the proof of~\cite[Theorem 3.4]{Langenbruch:2006ly} can be
adapted, with $M_{\alpha}^{A} \deq (A\Omega)^{\star}(2 (\alpha + 1)) /
2$, to establish Proposition~\ref{prop:ub}.  For the sake of
completeness, we produce our own version of this proof.  The third
lemma uses the previous one to bound the iterates of the Hermite
operator acting on $f$.

\begin{lemm}
    \label{lemm:Hermiteoperator}
    There exists a $K < \infty$ such that, for all weight functions
    $\omega$, there exists $C$, for all $f \in \Fspace_{\omega}$,
    \begin{equation*}
	\norm{L^2}{(x^2-D^2)^M f} \leq C {K}^M
	e^{\frac{\Omega^\star(4M+2)}{2}} \norm{\omega}{f}.
    \end{equation*}
\end{lemm}

\begin{proof}         
    Using Leibnitz' rule, we can expand $(x^2-D^2)^{M}$ as a sum of
    $2^M$ terms $C_{\mu,\nu} x^\nu D^\mu$ with $\mu+\nu \leq 2M$, and
    each $C_{\mu,\nu} \leq 4^M$.  Applying~(\ref{eq:numueff}) then
    yields
    \begin{equation}
	\norm{L^2}{(x^2-D^2)^{M} f} \leq C (8K)^M
	e^{\frac{\Omega^\star(4M+2)}{2}} \norm{\omega}{f} \nonumber
    \end{equation}
    which is the desired result.
\end{proof}

\begin{proof}[of Proposition~\ref{prop:ub}]    
    Let $\omega$ be fixed and let $f \in \Fspace_{\omega}$.  Recall
    that $h_{n}$ is also an eigenvector of the self-adjoint Hermite
    operator $x^2-D^2$, with eigenvalue $N \deq 2n+1$.  For $M \geq
    0$, we thus have
    \begin{equation*}
	\scal{(x^2-D^2)^M f}{h_n} = \scal{f}{(x^2-D^2)^M h_n} = N^M
	\scal{f}{h_n}
    \end{equation*}
    hence
    \begin{align*}
	\abs{\scal{f}{h_n}} &\leq N^{-M} \norm{L^2}{(x^2-D^2)^M f} \\
	\intertext{  using Lemma~\ref{lemm:Hermiteoperator} }
	&\leq C \paren{\frac{K}{N}}^{M}
	e^{\frac{\Omega^\star(4M+2)}{2}} \norm{\omega}{f} \\
	\intertext{ optimizing in $M$}
	&\leq C \paren{\frac{N}{K}}^{\frac12} e^{-\frac12 \sup_{M}
	\log\paren{\frac{N}{K}} (2 M + 1) - \Omega^\star(4 M + 2)}
	\norm{\omega}{f} \\
	\intertext{ since $\Omega$ is convex}
	&\leq C \paren{\frac{N}{K}}^{\frac12} e^{-\frac12
	\Omega\paren{\frac{1}{2} \log\paren{\frac{N}{K}} } }
	\norm{\omega}{f} \\
	\intertext{finally, using~\ref{it:minomega} and~\ref{it:croissomega}}
	& \leq C
	e^{-\frac{ \omega(\sqrt n)}{\Cclass}}
	\norm{\omega}{f}
    \end{align*}
    for some $\Cclass > 0$ that depends only on the constant
    in~\ref{it:croissomega}, thus only on $\bar\omega$.  In the
    previous reasoning we can then replace $\omega$ by $\Cclass
    \omega$ to obtain~(\ref{eq:ub}).
\end{proof}

\subsection{Lower bound}

To prove~(\ref{eq:lb}), we need a bound on the decrease of the Hermite
functions.  It is well known (see Szeg\"o's
reference~\cite[8.22.14]{Szego:1959uq} for instance) that
$\abs{h_{n}(x)} \leq 1$ for all $x\in\R, n \in \N$ (actually
$\abs{h_{n}(x)} \leq n^{-\frac{1}{12}}$).

\begin{lemm}
    \label{lemm:ineqdiff}
    Let $y$ and $z$ be two $C^2$ functions: $[x_0,+\infty) \to (0,
    +\infty)$ such that $y'(x) \to 0$, $z$ is bounded, satisfying the
    differential equations
    \begin{align*}
	y''(x) &= \phi(x) y(x) \\
	z''(x) &= \psi(x) z(x),
    \end{align*}
    with continuous $\phi(x) \leq \psi(x)$, and initial conditions
    $y(x_0) = z(x_0)$.  Then for all $x \geq x_0$, $z(x) \leq y(x)$.
\end{lemm}

\begin{proof} 
    Suppose that there exists $x_1 \geq x_0$ where $z(x_1) > y(x_1)$.
    Then for some $x_2 \in[x_0,x_1]$ we have $z'(x_2) > y'(x_2)$ and
    $z(x_2) \geq y(x_2)$.  Consequently, for all $x \geq x_2$, $z''(x)
    - y''(x) \geq 0$, and $z'(x)-y'(x) \geq z'(x_2) - y'(x_2)$.  When
    $x\to\infty$, $\liminf z'(x) \geq z'(x_2) - y'(x_2) > 0$, which
    contradicts the boundedness of $z$.
\end{proof}

\begin{lemm}
    \label{lemm:dh}
    For all $n \in \N$ and $\abs{x} \geq s \deq \sqrt{2n+1}$,
    $$\abs{h_n(x)} \leq {h_n(s)} e^{-\frac{(\abs{x}-s)^2}{2}} \leq
    e^{-\frac{(\abs{x}-s)^2}{2}}.$$
\end{lemm}

\begin{proof}
    By parity, we can assume $x \geq s$ (then $h_n(x) > 0$).  This
    implies that
    \begin{equation}
	\label{eq:ineqhn}
	(x-s)^2 -1 < x^2 - s^2.
    \end{equation}
    
    Recall that $h_n$ satisfies the differential equation $h_n'' =
    (x^2 - s^2) h_n$.  On the other hand, $y(x) \deq {h_n(s)}
    e^{-\frac{({x}-s)^2}{2}}$ satisfies $y'' = ((x-s)^2-1) y$.
    This, the obvious properties of $h_n$ and $y$,
    and~(\ref{eq:ineqhn}) together imply, by
    Lemma~\ref{lemm:ineqdiff}, that $h_n(x) \leq y(x)$.
\end{proof}

\begin{proof}[of Proposition~\ref{prop:lb}]
Since $\limsup \frac{\omega(t)}{t^2} < \frac12$, there exists 
$t_0 < \infty$ and  $\theta 
< 1$ such that, for all $t \geq t_0$, $\omega(t) \leq \theta^2 
\frac{t^2}{2}$. Let $n_{\theta}(x) \deq \frac{x^2 (1-\theta)^2 
    -1}{2}$.

We start with $\alpha \in \Lambda_{\Cclass \omega}$ ($\Cclass$ to be
determined later) and suppose that $\norm{\omega}{\alpha} \leq 1$,
which means that for all $n \geq 0$, $\abs{\alpha_{n}}\leq e^{-\Cclass
\omega(\sqrt n)}$.  Thus if $f = \sum_{n} \alpha_{n} h_{n}$ (this
series converging in $\Fspace$),
\begin{equation*}
    \abs{f(x)} \leq \underbrace{\sum_{0 \leq n \leq n_{\theta}(x) }
    e^{-\Cclass \omega(\sqrt n)} \abs{h_{n}(x)}}_{S_{1}} +
    \underbrace{\sum_{n > n_{\theta}(x)} e^{- \Cclass \omega(\sqrt n)}
    \abs{h_{n}(x)}}_{S_{2}}.
\end{equation*}

In the first sum, $\abs{x}-\sqrt{2n+1} \geq \theta \abs{x}$, so we can
use Lemma~\ref{lemm:dh} to bound $\abs{h_{n}(x)}$ by $e^{-\theta^2
\frac{x^2}{2}} \leq e^{-\omega(\abs{x})}$ as soon as $\abs x \geq
t_0$.  Because of~\ref{it:minomega}, the series $\sum_{n} e^{-\Cclass
\omega(\sqrt n)}$ converges, so $S_{1}$ is bounded by $C
e^{-\omega(\abs x)}$.

In the second sum, we simply bound $\abs{h_{n}(x)}$ by $1$ and use the
following on the tail of the sum.
\begin{align}
    \sum_{n \geq y} e^{-\Cclass \omega(\sqrt{n})} &\leq
    \int_{y-1}^{\infty} e^{-\Cclass\omega(\sqrt x)}d x \nonumber \\ &
    \leq \int_{\log(y-1)}^{\infty} e^{t-\Cclass\Omega(t/2)} dt \nonumber
    \\
     & \leq \frac{1}{\tfrac \Cclass 2 \Omega'(\tfrac{\log(y)}{2}) - 1}
     \int_{\log(y)}^{\infty} \paren{\tfrac \Cclass 2
     \Omega'(\tfrac{t}{2})-1} e^{t- \Cclass \Omega(\frac{t}{2})} dt
     \nonumber \\
     & \leq \frac{e^{\log(y)- \Cclass \omega(\sqrt{y})}}{\frac \Cclass
     2 \sqrt{y} \omega'(\sqrt{y}) -1} \nonumber \\
     & \leq \frac{\sqrt{y} e^{-\Cclass
     \omega(\sqrt{y})}}{\frac{\Cclass}{2} \omega'(\sqrt{y}) -
     \frac{1}{\sqrt{y}}} \nonumber \\
     &\leq C e^{-\frac{ \Cclass \omega(\sqrt y)}{ 2}}.
     \label{eq:mj}
\end{align}
We apply this with $y = n_{\theta}(x)$, $\omega(\sqrt y) \geq
\omega(\frac{1-\theta}{2} \abs{x}) \geq \frac{2}{\Cclass}
\omega(\abs{x})$ by~\ref{it:croissomega}, for some $\Cclass$ that
depends only on $\omegaclass$.  Finally we obtain $S_{2} \leq C
e^{-{\omega(x)}}$, and the proposition is proved.

\end{proof}

\paragraph{Remark.}
Proposition~\ref{prop:lb} can also be shown by adapting the proof
of~\cite[Theorem 3.4]{Langenbruch:2006ly}.

\section{Wigner distribution}
\label{sec:Wd}

We recall the definitions of $\Phi_{\rho}$ and $\Phi(f,f)$, already
given in the introduction as (\ref{eq:Wignerdistribution})
and~(\ref{eq:Wignertransform}).  If $\rho$ is a semi-definite positive
Hermitian operator: $L^2\to L^2$ that diagonalizes in an orthonormal
basis $\psi_{i}$, with eigenvalues $\prob_{i}$, then
\begin{equation}
    \nonumber \Phi_{\rho}(q,p) \deq \sum_{i} \frac{\prob_{i}}{2\pi}
    \int e^{ipx} \psi_{i}(q-\frac{x}{2})
    \overline{\psi_{i}(q+\frac{x}{2})} dx.
\end{equation}
If $\rho_{f}$ represents a pure state, in other words, if it is the
projector on the subspace generated by $f$, then we write
\begin{equation}
    \nonumber \Phi(f,f)(q,p) \deq \Phi_{\rho_{f}}(q,p) =
    \frac{1}{2\pi} \int e^{ipx} f(q-\frac{x}{2})
    \overline{f(q+\frac{x}{2})} dx.
\end{equation}

\subsection{Pure states}

The squared modulus of a function $f$ and its Fourier transform
$\Fourier{f}$ can be recovered from $\Phi(f,f)$ as marginal
distributions
\begin{equation}
    \int \Phi(f,f)(q,p) dp = \abs{f(q)}^2
    \label{eq:mq}
\end{equation}
and  
\begin{equation}
    \int \Phi(f,f)(q,p) dq = \bigabs{\Fourier{f}(p)}^2
    \label{eq:mp}.
\end{equation}

Let us define
\begin{equation}
    \norma{\omega}{\Phi} \deq \sup_{q,p} \abs{\Phi(q,p)}
    e^{2(\omega(\abs{q})+\omega(\abs{p}))} \nonumber
\end{equation}
and 
\begin{equation}
    \normb{\omega}{f} \deq
    \sqrt{\norma{\omega}{\Phi(f,f)}}.  \nonumber
\end{equation}
Clearly these are two norms, and using (\ref{eq:mq}) and
(\ref{eq:mp}) we get by integration, for all $q$ and
for all $p$,
\begin{equation}
    \abs{f(q)} e^{{\omega(\abs{q})}} \leq \sqrt{ \scriptstyle
    \int_{\R} e^{- 2 \omega(\abs{p})} dp } \normb{\omega}{f} \nonumber
\end{equation}
as well as
\begin{equation}
    \bigabs{\Fourier{f}(p)} e^{{\omega(\abs{p})}} \leq \sqrt{
    \scriptstyle \int_{\R} e^{- 2 \omega(\abs{q})} dq }
    \normb{\omega}{f} \nonumber
\end{equation}
the integrals being finite by~\ref{it:minomega}, hence finally
\begin{equation}
    \norm{ \omega }{f} \leq C \normb{\omega}{f}
    \label{eq:Wb}
\end{equation}
the constant $C$ depending only on $\omega$.  We now aim at
the corresponding lower bound.

\begin{prop}
    \label{prop:Wl}
    Assume~\ref{it:maxomegastrict}.  There exists $\Cclass < \infty$,
    depending only on $\omegaclass$, and $C < \infty$ (which may
    depend on $\omega$), such that for all $f \in \Fspace_{\Cclass
    \omega}$,
    \begin{equation}
	\normb{\omega}{f} \leq C \norm{\Cclass \omega}{f}.
	\label{eq:Wl}
    \end{equation}
\end{prop}

We prove actually a more general result, which gives (\ref{eq:Wl})
as a particular case.  Let
\begin{equation}
    \norm{\omega}{\rho} \deq \sup_{m,n} \abs{\rho_{m,n}}
    e^{\omega(\sqrt m) + \omega(\sqrt n)}.  \nonumber
\end{equation}

\begin{prop}
    \label{prop:Wlmixed}
    Assume~\ref{it:maxomegastrict}.  There exists $\Cclass < \infty$,
    depending only on $\omegaclass$, and $C < \infty$ (which may
    depend on $\omega$), such that for all $\rho$,
	\begin{equation}
	    \norma{\omega}{\Phi_{\rho}} \leq C \norm{\Cclass
	    \omega}{\rho}.
	    \label{eq:Wlmixed}
	\end{equation}
\end{prop}

The proof is based on radial bounds for the special Hermite
functions.  We recall that
\begin{equation}
    \Phi_{\rho}(q,p) = \sum_{m,n} \rho_{m,n} \Phi_{m,n}(q,p) \nonumber
\end{equation}
where $\Phi_{m,n}$, defined by (\ref{eq:defPhimn}), can also, as shown
in~\cite{Leonhardt:1997fk}, be expressed as follows: when $m \geq n$,
\begin{equation}
    \Phi_{m,n}(q,p) = \frac{(-1)^m}{\pi}
    \paren{\frac{n!}{m!}}^\frac12 e^{-\paren{q^2+p^2}} 
     \paren{\sqrt2(ip-q)}^{m-n}
    L_n^{m-n}\paren{2q^2+2p^2}.
    \label{eq:WmnLaguerre}
\end{equation}
Here $L_{n}^{\alpha} \deq (n!)^{-1} e^x x^{-\alpha} \frac{d^n}{dx^n}
(e^{-x}x^{n+\alpha})$ is the (non normalized) Laguerre polynomial of
degree $n$ and order $\alpha$.

If $m < n$, then by Hermitian symmetry $\Phi_{m,n}(q,p) =
\Phi_{n,m}(q,-p)$, which is equivalent to taking the canonical
generalization of Laguerre polynomials for $-n \leq \alpha < 0$:
\begin{math}
    L_n^\alpha(x) \deq \frac{(n+\alpha)!}{n!} (-x)^{-\alpha}
    L_{n+\alpha}^{-\alpha}(x).
\end{math}
But since $\rho$ is Hermitian, $\rho_{n,m} = \overline{\rho_{m,n}}$,
we only have to consider $m \geq n$ in the sums below.

The modulus of $\Phi_{m,n}$ is thus radial: writing $r \deq
\sqrt{q^2+p^2} $, we have
\begin{equation}
    \label{eq:deflmn}
    l_{m,n}(r) \deq \abs{\Phi_{m,n}(q,p)} =
    \frac{2^{\frac{m-n}2}}{\pi} \paren{\frac{n!}{m!}}^\frac12 e^{-r^2}
    r^{m-n} \abs{L_n^{m-n}(2r^2)}.
\end{equation}

What we need at this point is a bound on $l_{m,n}$ that is uniform on 
$m$ and $n$, in the same fashion as Lemma~\ref{lemm:dh}:

\begin{lemm}
    \label{lemm:boundLaguerre}
    There exists a constant $K$ such that, for all $m \geq n$ and $s
    \deq \sqrt{m+n+1}$, for all $r \geq 0$,
    \begin{equation}
	\label{eq:boundLaguerre}
	l_{m,n}(r) \leq K 
	\begin{cases}
	    1 &\text{ if } 0 \leq r \leq s \\
	    e^{-(r-s)^2} &\text{ if } r \geq s.
	\end{cases}
    \end{equation}
\end{lemm}

\begin{proof}
    When $r \leq s$, the result follows from the uniform bounds on
    Laguerre polynomials, for instance given by
    Krasikov~\cite{Krasikov:2007uq}:
    \begin{equation}
	\paren{L_{n}^\alpha (x)}^{2} e^{-x} x^{\alpha+1} \leq 1444
	n^{-\frac16} (n+\alpha+1)^{\frac12}
    \end{equation}
    to be used in (\ref{eq:deflmn}) with $x = 2r^2$ and $\alpha=m-n$.
    
    When $r \geq s$, $L_{n}^\alpha(2r^2)$ doesn't vanish and keeps the
    same sign as $L_{n}^\alpha(2s^2)$.  Now, as it can be seen
    from~\cite[5.1.2]{Szego:1959uq}, the function $z(r) \deq
    \sqrt{r}{l_{m,n}(r)}$ satisfies the differential equation
    $z''=(4(r^2-s^2)+\frac{\alpha^2-1/4}{r^2})z$.  On the other hand,
    $y(r) \deq \sqrt{s} {l_{m,n}(s)} e^{-(r-s)^2}$ satisfies $y'' =
    (4(r-s)^2-2) y$.  When $r \geq s$,
    \begin{equation}
	4(r-s)^2-2 < 4(r^2-s^2) + \frac{\alpha^2-1/4}{r^2}
    \end{equation}
    from which we conclude with Lemma~\ref{lemm:ineqdiff} that $z(r)
    \leq y(r)$.
\end{proof}

\begin{proof}[of Proposition~\ref{prop:Wl}]
    By (\ref{eq:rhomnforpurestates}) and
    Proposition~\ref{prop:ub}, if $\rho =
    \rho_{f}$ is a pure state, then for all $m,n$,
    \begin{equation}
	\abs{\rho_{m,n}}e^{\omega(\sqrt m)+\omega(\sqrt n)} =
	\abs{\scal{f}{h_{m}}} e^{\omega(\sqrt m)} \abs{\scal{f}{h_{n}}}
	e^{\omega(\sqrt n)} \leq C \norm{\Cclass \omega}{f}^2 \nonumber
    \end{equation}
    hence, applying Proposition~\ref{prop:Wlmixed},
    \begin{equation}
	\normb{\omega}{f}^2 = \abs{\Phi(f,f)}_{\omega} \leq
	C^{2} \norm{\Cclass^{2} \omega}{f}^2.  \nonumber
    \end{equation}
\end{proof}

\begin{proof}[of Proposition~\ref{prop:Wlmixed}]
    It is very similar to Proposition~\ref{prop:lb}.  There exists
    $t_{0} <\infty$ and $\theta < 1$ such that, for all $t \geq
    t_{0}$, $\omega(t) \leq \theta^2 \frac{t^2}{2}$.  Let
    $m_{\theta}(r) \deq {r^2 (1-\theta)^2 -1}$ and assume that
    $\norm{\Cclass \omega}{\rho} \leq 1$, $\Cclass$ to be chosen
    later.  For all $m,n$ we have $\abs{\rho_{m,n}} \leq
    e^{-\Cclass (\omega(\sqrt m) + \omega(\sqrt n))}$.  When ${m+n} \leq
    m_{\theta}(r)$, ${r}-s \geq \theta r$ and
    by~(\ref{eq:boundLaguerre}), this means that $l_{m,n}(r) \leq K
    e^{-\theta^2 {r^2}} \leq K e^{-2\omega(r)}$.  So
    \begin{equation}
	\label{eq:majsumagain1}
	\sum_{m+n \leq m_{\theta}(r)} \abs{\rho_{m,n}} l_{m,n}(r) \leq
	C e^{-2 \omega(r)}
    \end{equation}
    for $C \deq K \sum_{m,n} e^{ - \Cclass
    (\omega(\sqrt{m})+\omega(\sqrt{n})) }$.
    
    On the other hand, comparing the sum to the integral, we get
    similarly as in~(\ref{eq:mj})
    \begin{equation}
	\sum_{m+n \geq y} e^{ - \Cclass (\omega(\sqrt{m})+\omega(\sqrt{n}))
	} \leq C e^{-\frac{\Cclass}{\sqrt2} \omega(\sqrt y) }
	\nonumber
    \end{equation}
    using (\ref{eq:boundLaguerre}) again,
    \begin{align}
	\sum_{m+n \geq m_{\theta}(r) } \abs{\rho_{m,n} } l_{m,n}(r) &\leq 
	C e^{-\frac{\Cclass}{\sqrt{2}} 
	\omega(\sqrt{m_{\theta}(r)}) }
	\nonumber \\
	\label{eq:majsumagain2}
	&\leq C e^{- 2 \omega(r)}
    \end{align}
    if the constant $\Cclass$ is chosed large enough (depending only on
    $\omegaclass$).  Combining (\ref{eq:majsumagain1}) and
    (\ref{eq:majsumagain2}) yields the announced result.
    
\end{proof}

If the exact form of $\omega$ is known, such as $\omega(x) = x^\beta,
0 < \beta < 2$, a constant $\Cclass$ close to the optimal can easily
be obtained. An application to quantum statistics will be presented 
in a forthcoming paper.

\subsection{Mixed states}
\label{sec:mixedstates}

If we apply (\ref{eq:Wb}) to the right-hand side of
(\ref{eq:ub}), then use
(\ref{eq:rhomnforpurestates}), we obtain that the density matrix
coefficients of a pure state $f$ decrease as
\begin{equation}
    \abs{\rho_{m,n}} \leq C \normb{\Cclass \omega}{f}
    e^{-(\omega(\sqrt{m}) + \omega(\sqrt{n}))} \nonumber
\end{equation}
in other words, we have a converse to Proposition
\ref{prop:Wlmixed}, but only when $\rho = \rho_{f}$.  No such converse
can hold in the case of a general (mixed) state, as the following
example shows.

Following Butucea {\it et al.}~\cite{Butucea:2005uq}, let 
\begin{equation}
    \rho_{m,n} \deq \frac{(-1)^m \delta_{m,n}}{(m+1)(m+2)} =
    \begin{cases} 
	(-1)^m \int_{0}^{1} z^m (1-z) dz \text{ if } m = n \\
	0 \text{ else.}
    \end{cases}
    \label{eq:contrexample}
\end{equation}
We also recall an integral representation for Laguerre 
functions~\cite[5.4.1]{Szego:1959uq}, for $\alpha > -1$:
\begin{equation}
    L_{m}^{\alpha}(x) = \frac{e^{x} x^{\frac{\alpha}{2}}}{m!}
    \int_{0}^{\infty} e^{-t} t^{m+\frac{\alpha}{2}}
    J_{\alpha}(2\sqrt{tx}) dt.  \nonumber
\end{equation}
Then, using (\ref{eq:WmnLaguerre}) once again with $\alpha = m - n =
0$ and $x \deq 2({q^2+p^2})$,
\begin{align*}
    \Phi_{\rho}(q,p) &= \frac{e^{-\frac x 2}}{\pi} \sum_{m}
    \int_{0}^{1} z^m (1-z) dz L_{m}^{0}(x) \\
    &= \frac{e^{\frac x 2}}{\pi} \int_{0}^{1} \int_{0}^{\infty}
    \sum_{m} z^m (1-z) \frac{t^m}{m!} e^{-t} J_{0}(2\sqrt{t x}) dt dz
    \\
    &= \frac{e^{\frac x 2}}{\pi} \int_{0}^{1} (1-z) \int_{0}^{\infty}
    e^{(z-1)t} J_{0}(2\sqrt{tx}) dt dz \\
    &= \frac{e^{\frac x 2}}{\pi} \int_{0}^{1} e^{\frac{x}{z-1}} dz \\
    &\leq \frac{e^{-\frac x 2}}{\pi}
\end{align*}
We see that, although $\rho_{m,n}$ has algebraic decrease along the
diagonal, $\Phi_{\rho}(q,p)$ decreases like $e^{-(q^2+p^2)}$.
Cancellations occur because of the alternating signs in $\rho_{m,m}$.

\subsection{Ambiguity function}
\label{sec:ambiguityfunction}

The \emph{ambiguity function} $A(f,f)$ of a signal $f \in L^2$ is
simply the Fourier transform in $(q,p)$ of its Wigner transform.
Knowing (\ref{eq:FourierspecialHermite}),
\begin{equation}
    A(f,f)(\vartheta, \varpi) \deq \Fourier{\Phi(f,f)}(\vartheta,
    \varpi) = \sum_{m,n} \frac{(-i)^{m+n}}{2} \scal{f}{h_{m}}
    \overline{\scal{f}{h_{n}}} \Phi_{m,n}\paren{\tfrac{\vartheta}{2},
    \tfrac{\varpi}{2}} \nonumber
\end{equation}

The proof of Proposition~\ref{prop:Wl} works \textit{mutatis
mutandis} to obtain
\begin{equation}
    \sqrt{\norma{\omega}{A(f,f)}} \leq C \norm{\Cclass \omega}{f}
    \nonumber
\end{equation}
the constant $C$ depending on $\omega$ and $\Cclass$ depending on
$\omegaclass$.

Naturally, the equivalent of (\ref{eq:Wb}) holds also:
\begin{equation}
    \norm{\omega}{f} \leq C \sqrt{\norma{\Cclass \omega}{A(f,f)}}
    \nonumber
\end{equation}
in particular $\Phi(f,f)$ and $A(f,f)$ are tamely equivalent in their
respective $\omega$ norms.  However, as the previous example
(\ref{eq:contrexample}) shows, this is not necessarily true for mixed
states.

\begin{acknowledgements}
    The author would like to thank the anonymous referee for his 
    valuable comments and suggestions.
\end{acknowledgements}

\affiliationone{
Jean-Marie Aubry \\
Universit\'e Paris-Est \\
LAMA -- UMR CNRS 8050 \\
94010 CR\'ETEIL\\
FRANCE \email{jmaubry@math.cnrs.fr} }

\end{document}